\documentclass[11pt]{amsart}
\usepackage[margin=1.2in]{geometry}
\usepackage{graphicx}
\usepackage{amssymb}
\usepackage{amsthm}
\usepackage{epstopdf}
\usepackage[colorlinks=true, pdfstartview=FitV, linkcolor=blue, citecolor=blue, urlcolor=blue]{hyperref}
 \usepackage{tikz}
\usepackage[listings]{tcolorbox}
\usepackage{caption}
\usepackage{subcaption}
\usepackage{mathrsfs}
\usepackage{enumitem}
\usepackage{soul}
\usepackage{verbatim}
\usepackage{faktor}
\usepackage{bbm}
\usepackage{ytableau}
\usepackage[nameinlink]{cleveref}
\usepackage{yhmath}
\usepackage{biblatex}
\usepackage{color}

\newcommand{\s}{\sigma}

\newcommand{\N}{{\mathbb N}}
\renewcommand{\P}{{\mathbb P}}

\newcommand{\ba}{{\mathbf a}}

\newcommand{\bx}{{\mathbf x}}

\def\Tor{\operatorname{Tor}}

\def\Gr{\operatorname{Gr}}

\renewcommand{\S}{\mathfrak{S}}

\theoremstyle{plain} 
\newtheorem{thm}{Theorem}[section]

\newtheorem*{introthm*}{Theorem}
\newtheorem*{introcor*}{Corollary}

\newtheorem{lem}[thm]{Lemma}
\newtheorem{prop}[thm]{Proposition}

\newtheorem*{oprobl*}{Open Problem}

\newtheorem{conj}[thm]{Conjecture}

\theoremstyle{definition}
\newtheorem{defn}[thm]{Definition}

\newtheorem{constr}[thm]{Construction}

\theoremstyle{remark}
\newtheorem*{rem}{Remark}

\AddToHook{env/defn/begin}{\crefalias{thm}{defn}}
\AddToHook{env/lem/begin}{\crefalias{thm}{lem}}
\AddToHook{env/cor/begin}{\crefalias{thm}{cor}}
\AddToHook{env/prop/begin}{\crefalias{thm}{prop}}
\AddToHook{env/conj/begin}{\crefalias{thm}{conj}}
\AddToHook{env/constr/begin}{\crefalias{thm}{constr}}
\AddToHook{env/section/begin}{\crefalias{thm}{section}}

\numberwithin{equation}{section}

\title{The Minimum Number of Generators of Symmetric Ideals}
\author{Noah Walker}

\begin{document}

\begin{abstract}
We study equigenerated symmetric ideals in polynomial rings and the minimum number of polynomials required to generate them up to permutations of the variables.
We give a representation-theoretic formula for this number and determine a sharp threshold on the number of variables needed for a general symmetric ideal to attain the largest possible dimension in its generating degree.
The threshold is governed by partial sums of integer partition numbers, which appear in the OEIS as sequence A000070. 
We also construct explicit extremal symmetric ideals for every admissible number of generators. 
As an application, the principal case gives the sharp stable range for the theorem of Harada–Seceleanu–\c{S}ega on general principal symmetric ideals in \cite{HSS}.
\end{abstract}

\maketitle

\setcounter{tocdepth}{1}
\tableofcontents

\section{Introduction}

The study of symmetric ideals in commutative algebra explores a special class of ideals that remain invariant under permutations of variables, typically within polynomial rings.
These ideals arise naturally in contexts where underlying algebraic structures exhibit symmetry, such as in invariant theory, algebraic statistics, and representation theory \cite{FiniteGen, EquivariantGrobner, FIMods, InvariantHilbertSchemes, FinitenessInvariantChains, AsymptoticSymmetric}.

Let $k$ be a field and let $R=k[x_1,\dots,x_n]$.
Let $\S_n$ be the symmetric group on $n$ letters.
Then $\S_n$ acts on $R$ by permuting the variables.
That is, for $\s\in \S_n$,
$$\s\cdot f(x_1,\dots,x_n)=f(x_{\s(1)},\dots,x_{\s(n)}).$$
If an ideal $I$ of $R$ is stable under this action, then we call $I$ a {\bf symmetric ideal}.
For $f_1,\dots,f_r\in R$, the symmetric ideal generated by $f_1,\dots,f_r$ is defined to be
$$(f_1,\dots,f_r)_{\S_n}=(\s\cdot f_i : \s\in \S_n, 1\leq i\leq r).$$
The case $r=1$ gives the class of \emph{principal symmetric ideals}.
In this paper, we study the minimum number $r$ such that a symmetric ideal can be generated by $r$ elements.
We approach this via a representation-theoretic framework, studying the minimum number of $k\S_n$-module generators of the minimum degree component of a symmetric ideal.
In particular, we obtain the following formula for this number.

\begin{introthm*}[\Cref{r gen characterization}]
	Let $V\cong\bigoplus_{i=1}^t n_iV_i$ be a finitely generated $k\S_n$-module where $V_1,\dots,V_t$ are pairwise non-isomorphic simple $k\S_n$-modules.
	Then the smallest number of generators for $V$ as a $k\S_n$-module is
    $$r=\left\lceil\max\left\{\frac{n_i}{\dim_k V_i}\mid 1\leq i\leq t\right\}\right\rceil.$$
\end{introthm*}

Specializing to $r=1$ gives a recognition theorem for principal symmetric ideals, answering a question raised in \cite{PolymathJr}.
The result yields an algorithm for computing the minimal number of orbit generators of a symmetric ideal from the canonical decomposition of the associated $k\S_n$-module. 

For general symmetric ideals generated in a fixed degree $d$ with a fixed number of orbit generators $r$, their behavior is often extremal \cite{HSS, SS}.
See \Cref{general symmetric ideal} for a description of these ideals.
In \Cref{general r gen is max}, we prove that such ideals attain the largest possible dimension in their generating degree given the number of orbit generators.
This motivates us to study submodules of a fixed $k\S_n$-module which are extremal with respect to a fixed number of generators $r$.
We call such modules {\em maximal $r$-generated submodules} (see \Cref{def MrGS}) and determine their structure.

\begin{introthm*}[\Cref{maximal r gen sub decomp}]
	Let $V\cong\bigoplus_{i=1}^t n_iV_i$ be a finitely generated $k\S_n$-module where $V_1,\dots,V_t$ are pairwise non-isomorphic simple $k\S_n$-modules.
	For a fixed $r$, every maximal $r$-generated submodule of $V$ is isomorphic to 
	$$\bigoplus_{i=1}^t\min\{n_i,r\cdot\dim_k V_i\}V_i.$$
	In particular, maximal $r$-generated submodules are unique up to isomorphism.
\end{introthm*}

Applying this result to the homogeneous component $R_d$ naturally leads to sums of Kostka numbers.
To analyze these sums, we prove a new inequality (\Cref{main inequality}) comparing their growth with the dimensions of Specht modules.
This inequality identifies the isotypic component of $R_d$ which requires the largest number of generators as a $k\S_n$-module.
In \Cref{kostka partition sum}, we also give a combinatorial proof identifying some of these sums of Kostka numbers as the OEIS sequence \href{https://oeis.org/A000070}{A000070}, which is given by $\sum_{i=0}^{d-1}P(i)$ where $P(i)$ is the number of integer partitions of $i$.
Besides providing a key identity for our main theorem, this identity gives the sequence a representation-theoretic interpretation, which seems not to have been observed previously.
These combinatorial tools allow us to compute the dimension of maximal $r$-generated submodules of $R_d$ in the following theorem.
\begin{introthm*}[\Cref{main theorem}, \Cref{general r gen is max}]
	Fix $n>d\ge 2$ and $1\le r\le P(d)$, where $P(d)$ denotes the number of partitions of $d$.
	Let $M$ be a maximal $r$-generated submodule of $R_d$.
	Then $\dim_k M\le \dim_k R_d - (P(d)-r)$ with equality if and only if
	\begin{equation}\label{intro bound}
	n \ge 1+\frac{\sum_{i=0}^{d-1}P(i)}{r}.
	\end{equation}
	Consequently, a general $(r,d)$-symmetric ideal $I$ satisfies $\dim_k I_d=\dim_k R_d-(P(d)-r)$ if and only if \eqref{intro bound} holds.
\end{introthm*}

Finally, we construct explicit examples of extremal $r$-generated symmetric ideals achieving the maximal possible dimension when \eqref{intro bound} holds.

\begin{introthm*}[\Cref{J props}]
	Fix $n>d\ge 2$ and $1\le r\le P(d)$.
	Then there exists a symmetric ideal $J$ generated in degree $d$ with $\dim_k J_d=\dim_k R_d-(P(d)-r)$.
	The ideal $J$ is generated by $r$ elements up to symmetry if and only if \eqref{intro bound} holds.
\end{introthm*}

As an application, we use these examples to provide an optimal bound on the number of variables required for the theorem of Harada, Seceleanu, and \c{S}ega in \cite{HSS} to hold.
We note that \cite{SS} proves an analogous theorem and bound for general $(r,d)$-symmetric ideals using the results of this paper.

\begin{introcor*}[\Cref{8.4 strengthening}]
	Fix $n>d\ge2$.
    Let $I$ be a general principal symmetric ideal.
    Then the Hilbert function, Betti table, and $\S_n$-equivariant minimal free resolution of $R/I$ given by \cite[Theorem 8.4]{HSS} occur if and only if \eqref{intro bound} holds with $r=1$.
\end{introcor*}

The paper is organized as follows.
\Cref{characterization section} develops the characterization of $r$-generated symmetric ideals and introduces maximal $r$-generated submodules.
\Cref{symmetric rep theory} recalls the necessary representation theory of the symmetric group and
establishes the combinatorial identity involving Kostka numbers.
\Cref{kostka inequality section} proves the inequality needed to determine the structure of maximal $r$-generated submodules of $R_d$.
Finally, \Cref{gsi section} applies these results to general principal symmetric ideals and establishes the sharp bound on the number of variables required in the theorem of Harada, Seceleanu, and \c{S}ega.

Throughout the paper, we assume that $k$ is an infinite field with $\operatorname{char}(k)=0$ or $\operatorname{char}(k)>n$, so that $k\S_n$ is semisimple by Maschke's theorem.
Unless otherwise stated, all ideals are homogeneous and generated in a single degree.

\section{Characterization of $r$-generated symmetric ideals}\label{characterization section}

\begin{defn}\label{psi defn}
    Let $I$ be a symmetric ideal of $R$.
    We say that $I$ is {\bf $r$-generated} if there exist $r$ homogeneous polynomials $f_1,\dots,f_r\in R$ of the same degree such that $I=(f_1,\dots,f_r)_{\S_n}$.
	When $r=1$, we call $I$ a {\bf principal symmetric ideal}.
\end{defn}
\begin{defn}\label{r-gen defn}
	A $k\S_n$-module $M$ is {\bf $r$-generated} if there exist $f_1,\dots,f_r\in M$ such that $M=\langle f_1,\dots,f_r\rangle_{k\S_n}$, where $\langle f_1,\dots,f_r\rangle_{k\S_n}$ is the $k\S_n$-module generated by $f_1,\dots,f_r$.
\end{defn}

A symmetric ideal $I$ generated in degree $d$ is $r$-generated if and only if $I_d$ is $r$-generated as a $k\S_n$-module.

\begin{defn}
    Let $V_1,\dots,V_t$ be a complete set of representatives of the isomorphism classes of simple $k\S_n$-modules.
    Let $V$ be a finitely generated $k\S_n$-module.
    We call $V\cong\bigoplus_{i=1}^t n_i V_i$ the {\bf canonical decomposition} of $V$, where $n_i V_i$ denotes the direct sum of $n_i$ copies of $V_i$.
\end{defn}

\begin{thm}\label{r gen characterization}
    Let $V$ be a finitely generated $k\S_n$-module with canonical decomposition $V\cong\bigoplus_{i=1}^t n_i V_i$.
    Then $V$ is an $r$-generated $k\S_n$-module if and only if $n_i\leq r\cdot\dim_k V_i$ for each $i$.
    In particular, the smallest $r$ such that $V$ is $r$-generated is 
    $$r=\left\lceil\max\left\{\frac{n_i}{\dim_k V_i}\mid 1\leq i\leq t\right\}\right\rceil.$$
\end{thm}
\begin{proof}
    A $k\S_n$-module is $r$-generated precisely when it is a quotient of $r(k\S_n)$.
    In our setting, $\operatorname{char}(k)=0$ or $\operatorname{char}(k)>n$, so by the semisimplicity of $k\S_n$, this is equivalent to being isomorphic to a submodule of $r(k\S_n)$.
    The canonical decomposition of $r(k\S_n)$ is
    $$\bigoplus_{i=1}^t\ (r\cdot \dim_k V_i)V_i$$
	giving the desired result.
\end{proof}

Since the canonical decomposition of $I_d$ can be computed algorithmically, this characterization provides a method to compute the minimal $r$ such that a symmetric ideal is $r$-generated (see Macaulay2 code in \cite{GithubRepo} for an implementation).

\begin{defn}\label{def MrGS}
    A {\bf maximal $r$-generated submodule} of a module $V$ is an $r$-generated submodule of $V$ that is not contained in any other $r$-generated submodule of $V$.
\end{defn}
Throughout, all maximal $r$-generated submodules are as $k\S_n$-submodules.
In this case, the following proposition gives a complete description of their structure.

\begin{prop}\label{maximal r gen sub decomp}
    Let $V$ be a finitely generated $k\S_n$-module with canonical decomposition $V\cong\bigoplus_{i=1}^t n_i V_i$.
    Then every maximal $r$-generated submodule of $V$ is isomorphic to $$\bigoplus_{i=1}^t \min\{n_i,r\cdot\dim_k V_i\} V_i.$$
    In particular, maximal $r$-generated submodules of $V$ exist and are unique up to isomorphism.
\end{prop}
\begin{proof}
	Since $\bigoplus_{i=1}^t \min\{n_i,r\cdot\dim_k V_i\} V_i$ is $r$-generated by \Cref{r gen characterization}, it suffices to show that any $r$-generated submodule of $V$ is isomorphic to a submodule of $\bigoplus_{i=1}^t \min\{n_i,r\cdot\dim_k V_i\} V_i$.

    Let $M$ be an $r$-generated submodule of $V$ and let $M\cong\bigoplus_{i=1}^t m_i V_i$ be the canonical decomposition of $M$.
    Since $M$ is a submodule of $V$, $m_i\leq n_i$.
	By \Cref{r gen characterization}, $m_i\leq r\cdot\dim_k V_i$, so $m_i\leq\min\{n_i,r\cdot\dim_k V_i\}$.
    Thus, $$M\cong\bigoplus_{i=1}^t m_iV_i\subseteq\bigoplus_{i=1}^t\min\{n_i, r\cdot\dim_k V_i\}V_i. \qedhere $$
\end{proof}

\section{Maximal $r$-generated submodules of $R_d$ and Kostka numbers}\label{symmetric rep theory}

A {\bf partition} of $n$ is a nonincreasing ordered tuple of nonnegative integers $\lambda=(\lambda_1,\lambda_2,\dots,\lambda_m)$ such that $\sum_{i=1}^m\lambda_i=n$.
We call $\lambda_1,\lambda_2,\dots,\lambda_m$ the {\bf parts} of $\lambda$ and denote by $\#\lambda$ the number of parts of $\lambda$.
The notation $\lambda\vdash n$ means that $\lambda$ is a partition of $n$.
Generally, the parts of a partition are positive, but 
a partition may be padded with 0's to have the desired number of parts.
For $\lambda,\mu\vdash n$, $\lambda$ {\bf dominates} $\mu$, denoted $\lambda\unrhd\mu$, if for all $j$, $\sum_{i=1}^j \lambda_i\geq\sum_{i=1}^j \mu_i.$
The {\bf diagram} of $\lambda$ is an upper-left justified array of boxes such that row $i$ has length $\lambda_i$.
For two partitions $\nu,\lambda$, we write $\nu\subseteq\lambda$ if $\nu_i\le\lambda_i$ for all $i$.
In this case, the {\bf skew diagram} $\lambda/\nu$ consists of the boxes in $\lambda$, but not in $\nu$.
A skew diagram is a {\bf horizontal strip} if each column contains at most one box.

Recall that the {\bf Specht modules}, denoted by $S^\lambda$ for $\lambda\vdash n$, are precisely the simple $k\S_n$-modules.
The vector space dimension of $S^\lambda$, which we denote $f^\lambda$, can be interpreted combinatorially as the number of standard Young tableaux of shape $\lambda$.
A {\bf standard Young tableau} (SYT) of shape $\lambda$ is a filling of the boxes of the diagram of $\lambda$ with the numbers $1,\dots,n$ with no repeats so that the rows and columns are strictly increasing to the right and down, respectively.
A {\bf semistandard Young tableau} (SSYT) of shape $\lambda$ and content $\mu$, where $\lambda,\mu\vdash n$, is a filling of the diagram of $\lambda$ with $\mu_1$ 1's, $\mu_2$ 2's, etc. such that the rows weakly increase to the right and the columns strictly increase going down.
The number of SSYT of shape $\lambda$ and content $\mu$ is the {\bf Kostka number} of $\lambda$ and $\mu$, denoted $K_{\lambda,\mu}$.
An important class of modules in the representation theory of $\S_n$ is the class of permutation modules. 
For $\mu\vdash n$, its associated permutation module is $M^\mu:=\bigoplus_{\lambda\vdash n}K_{\lambda,\mu}S^\lambda$.
See \cite{James} for more on the above definitions.

\begin{defn}
    Let $\ba=(a_1^{m_1},a_2^{m_2},\dots,a_\ell^{m_\ell})\vdash d$ with $n$ parts, where the exponent denotes repetition of a part.
    The partition of $n$ obtained by reordering the vector $(m_1,m_2,\dots,m_\ell)$ so that the entries are nonincreasing is called the {\bf metatype} of $\ba$ and is denoted by $m(\ba)$.
\end{defn}
For a partition $\ba=(a_1,\dots,a_n)$ with $n$ parts, we use $\bx^\ba$ to denote the monomial $x_1^{a_1}\cdots x_n^{a_n}$.
It is a standard result that $M^{m(\ba)}$ is isomorphic to $\langle \bx^\ba\rangle_{k\S_n}$ as a $k\S_n$-module.
Combined with \Cref{maximal r gen sub decomp}, this gives the following result.
\begin{prop}\label{decomp of R_d}
    The canonical decomposition of $R_d$ as a $k\S_n$-module is
    $$R_d\cong\bigoplus_{\lambda\vdash n}\left(\sum_{\ba\vdash d,\#\ba=n} K_{\lambda,m(\ba)}\right)S^\lambda$$
    and any maximal $r$-generated submodule of $R_d$ is isomorphic to 
    \begin{equation}\label{kostka min}
    \bigoplus_{\lambda\vdash n}\min\left\{\sum_{\ba\vdash d,\#\ba=n}K_{\lambda,m(\ba)},r f^\lambda\right\}S^\lambda.
    \end{equation}
\end{prop}

The sums in \eqref{kostka min} will arise frequently, so we introduce the notation
$$F^\lambda_d:=\sum_{\ba\vdash d,\#\ba=n}K_{\lambda,m(\ba)}.$$
The number $F^\lambda_d$ can be interpreted combinatorially as the number of column-strict plane partitions of shape $\lambda$ and sum $d$.
We study these sums further in \Cref{kostka inequality section}, where we compute the minimum in \eqref{kostka min}.
We record here some elementary results about Kostka numbers and the dimensions $f^\lambda$.
Properties (1)-(4) are clear from counting tableaux.
Property (5) is well-known (see e.g. \cite[Sec 1.7]{MacDonald}).

\begin{prop}\label{kostka props}
    Let $\lambda,\mu\vdash n$.
    \begin{enumerate}
        \item $K_{\lambda,(1^n)}=f^\lambda$.
        \item $K_{(n),\mu}=1$.
        \item $f^{(n)}=1$.
        \item $f^{(n-1,1)}=n-1$.
        \item $K_{\lambda,\mu}>0$ if and only if $\lambda\unrhd\mu$.
    \end{enumerate}
\end{prop}

\begin{prop}\label{n-1 kostka}
    Let $n\geq 2$ and $\mu\vdash n$.
    If the parts of $\mu$ are all positive, then $K_{(n-1,1),\mu}=\#\mu-1$.
	In particular, $K_{(n-1,1),m(\ba)}=\#m(\ba)-1$ for any partition $\ba$.
\end{prop}
\begin{proof}
    Consider the diagram for $(n-1,1)$:
    
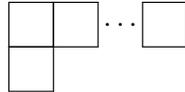
\begin{ytableau}
        \phantom{.}&\phantom{.}& \none[\cdots] & \phantom{.} \\
        \phantom{.}
    \end{ytableau}\ .
    The tableau is completely determined by the choice of the bottom row entry since the remaining integers must be placed in increasing order along the top row.
    Furthermore, any of the numbers $2,3,\dots,\#\mu$ can be chosen for the bottom row.
    Hence, $K_{(n-1,1),\mu}=\#\mu-1$.
\end{proof}

Recall that $P(i)$ denotes the number of partitions of $i$.
The sequence $\sum_{i=0}^{d-1}P(i)$ from the following proposition appears in the OEIS as A000070 and also appears in our bound in our main theorem.
There are many problems in which the sequence arises, primarily in partition theory \cite{OEIS}.
Its appearance in commutative algebra is perhaps surprising, but the following result shows that it appears as a sum of Kostka numbers of the form from \Cref{decomp of R_d}.

\begin{prop}\label{kostka partition sum}
    If $n>d$, then
    $$F^{(n-1,1)}_d=\sum_{i=0}^{d-1}P(i).$$
\end{prop}
\begin{proof}
	We show that both quantities count the sum over all partitions of $d$ of the number of distinct positive parts in each partition.
	(For example, for $d=3$, the partitions of $d$ are $(3)$, $(2,1)$, and $(1,1,1)$, so the sum is $1+2+1=4$.)
	This interpretation of the right hand side is standard (see e.g. \cite{OEIS} or \cite[Exercise 3.4.36]{West}).

	Consider the left hand side.
    By \Cref{n-1 kostka}, $K_{(n-1,1),m(\ba)}=\#m(\ba)-1$.
    Each part in $m(\ba)$ corresponds to a distinct part in $\ba$, so $\#m(\ba)$ is the number of distinct parts in $\ba$.
    Since $n>d$, $\ba$ always contains at least one 0, so $\#m(\ba)-1$ is the number of distinct positive parts in $\ba$.
    Thus, $F^{(n-1,1)}_d=\sum_{\ba\vdash d,\#\ba=n}K_{(n-1,1),m(\ba)}$ is the sum of the numbers of distinct positive parts in each partition of $d$.
\end{proof}

\section{Dimension of maximal $r$-generated submodules}\label{kostka inequality section}

The goal of this section is to show that when $n$ is large, the minimum in \eqref{kostka min} evaluates to the sum of Kostka numbers (except when $\lambda=(n)$).
This implies that $R_d$ can be decomposed into a maximal $r$-generated submodule and some copies of the trivial representation.
In \Cref{gsi section}, we identify these components explicitly.
We evaluate the minimum by showing that the ratio $\frac{F^\lambda_d}{f^\lambda}$ is largest when $\lambda=(n-1,1)$ and that the ratio for $\lambda=(n-1,1)$ is bounded above by $r$ when $n$ is sufficiently large.
To prove 
$$\frac{F^\lambda_d}{f^\lambda}\le\frac{F^{(n-1,1)}_d}{f^{(n-1,1)}}$$
it suffices to show 
\begin{equation}\label{kostka term inequality}
(n-1)K_{\lambda,\mu}\leq (\#\mu-1)f^\lambda.
\end{equation}
Unfortunately, \eqref{kostka term inequality} does not hold for $\lambda=\mu=(2,2)$, so we have to do some technical case analysis.
For the main case, we argue by induction using the following well-known recursion.
This is the usual recursion for Kostka numbers that corresponds to removing the boxes with the largest value.
\begin{lem}\label{kostka recursion}
    Let $\lambda,\mu\vdash n$ where $\mu=(\mu_1,\dots,\mu_\ell)$.
    Then
    $$K_{\lambda,\mu}=\sum_\nu K_{\nu,\mu'},$$
    where $\mu'=(\mu_1,\dots,\mu_{\ell-1})$ and the sum is taken over all partitions $\nu\subseteq\lambda$ such that $\lambda/\nu$ is a horizontal strip of size $\mu_\ell$.
\end{lem}

\begin{lem}\label{dimension inequality}
    Let $\lambda\vdash n$ and $s\in\N$.
    Then 
    $$\sum_\nu f^\nu\leq f^\lambda$$
    where the sum is taken over all partitions $\nu\subseteq\lambda$ such that $\lambda/\nu$ is a horizontal strip of size $s$.
\end{lem}
\begin{proof}
	Denote the set of SYT of shape $\lambda$ by $SYT(\lambda)$.
	Given $t\in SYT(\nu)$ for some $\nu$ satisfying the conditions, we can define a SYT $t'$ of shape $\lambda$ by filling each box in $\nu$ with its value from $t$ and filling the new boxes with the remaining $s$ values so that they increase along rows.
This defines an injection from $\bigcup_\nu SYT(\nu)$ to $SYT(\lambda)$, so the inequality holds.
\end{proof}

\begin{lem}\label{kostka inequality}
	Let $\lambda,\mu\vdash n$ with $\lambda\not=(n)$ and either $\lambda\ne(2,2)$ or $\mu\ne(2,2)$.
	Then
	$$(n-1)K_{\lambda,\mu}\le (\#\mu-1)f^\lambda.$$
\end{lem}
\begin{proof}
	Let $\ell=\#\mu$ and $m=\#\lambda$.
    
    {\em Case 1:} If any of the following conditions holds: $\lambda\not\!\unrhd \ \mu$, $\lambda=(n-1,1)$, $\mu=(n)$, $\mu=(1^n)$, or $\lambda=(1^n)$, then the result is immediate from \Cref{kostka props} and \Cref{n-1 kostka}.

    {\em Case 2:} Suppose that $n\leq 4$.
    The few cases not covered by Case 1 can be checked by hand.

    {\em Case 3:} Suppose that $\mu=(2,2,2)$ or $\mu=(2,2,1)$.
    Again, we can manually check that the inequality holds for all $\lambda$.
    The details are given in \Cref{kostka inequality appendix}.

    {\em Case 4:} Suppose that $m=2$ and $\lambda_2\leq\mu_\ell$.
	Then $\lambda=(n-p,p)$ for some $p\geq 2$.
    The hook lengths of $\lambda$ follow the straightforward pattern illustrated here.

    \ytableausetup{boxsize=3em}
    \begin{center}
    \begin{ytableau}
        \scriptstyle n-p+1 & n-p & \none[\cdots] & \scriptstyle n-2p+3 & \scriptstyle n-2p+2 & n-2p & \none[\cdots] & 2 & 1 \\
        p & p-1 & \none[\cdots] &  2 & 1\\
    \end{ytableau}
    \end{center}
    \normalsize
    
    Thus, by the hook length formula,
    $$f^\lambda=\frac{n!}{\prod h_\lambda(i,j)}=\frac{n!}{p!(n-2p)!\frac{(n-p+1)!}{(n-2p+1)!}}=\frac{n!(n-2p+1)!}{p!(n-p+1)!(n-2p)!}=\frac{n!(n-2p+1)}{p!(n-p+1)!}.$$

    To determine $K_{\lambda,\mu}$, note that these semistandard $\lambda$-tableaux are completely determined by the content in the bottom row of $\lambda$.
    There are $p$ boxes in the bottom row of $\lambda$ and $\ell-1$ values.
    The boxes are indistinguishable (since the bottom row must be weakly increasing) and the values are distinguishable.
    Any assignment of values to boxes is valid since $\lambda_2\leq\mu_\ell\leq\mu_1$, so $K_{\lambda,\mu}$ is counting the number of ways to assign $p$ indistinguishable objects to $\ell-1$ distinguishable categories.
    Thus, $K_{\lambda,\mu}=\binom{p+\ell-2}{\ell-2}$.

    The inequality that we are trying to show is thus
    $$\frac{(n-1)(p+\ell-2)!}{p!(\ell-2)!}\leq\frac{(\ell-1)n!(n-2p+1)}{p!(n-p+1)!}$$
    where we may assume $\ell\geq 2, p\geq 2, n\ge 5$ by the previous cases.
    Also,
    $n\geq \ell\cdot\mu_\ell\geq \ell\cdot\lambda_2=\ell\cdot p.$
    With these properties ($\ell,p\geq2$ and $n\geq5$ and $n\geq \ell\cdot p$), the inequality now follows from elementary estimates and calculations, the details of which are given in \Cref{inequality proof}.

    {\em Case 5:} Suppose that either $m=2$ and $\lambda_2>\mu_\ell$ or $m>2$.
    Also assume that $\mu\not=(2,2,2),(2,2,1)$.
    	We proceed by induction on $n$.
    We first show that for each $\nu$ in the sum in \Cref{kostka recursion}, $\nu\ne(n-\mu_\ell)$.
    	
    	If $\nu$ is contained in the first row of $\lambda$, then every other row of $\lambda$ is part of $\lambda/\nu$.
    	If $m>2$, then at least two rows of $\lambda$ are part of $\lambda/\nu$, so $\lambda/\nu$ is not a horizontal strip.
    	If $m=2$ and $\lambda_2>\mu_\ell$, then $\lambda/\nu$ has size $\ge \lambda_2>\mu_\ell$.
    	Hence, $\nu$ is not contained in a single row, so $\nu\ne(n-\mu_\ell)$.

    Note that if $\mu'=(2,2)$, then $\mu=(2,2,2)$ or $\mu=(2,2,1)$.
    However, by assumption, $\mu$ is neither of those, so $\mu'\not=(2,2)$.
    Thus, by the induction hypothesis,
    \begin{equation}\label{eqn:IH}
        (n-\mu_\ell-1)K_{\nu,\mu'}\leq (\#\mu'-1)f^\nu=(\ell-2)f^\nu.
    \end{equation}
    Since $\mu_i\geq\mu_\ell$ for all $i$, $n\geq \mu_\ell\cdot \ell$.
    Hence,
    $$\mu_\ell(\ell-2)=\mu_\ell\cdot\ell-2\mu_\ell\leq n-\mu_\ell-1,$$
    By Case 1, $\ell\ge2$.
    If $\ell\ge 3$,
    \begin{equation}\label{eqn:secondTerm}
        \mu_\ell K_{\nu,\mu'}\leq \frac{n-\mu_\ell-1}{\ell-2}K_{\nu,\mu'}\leq f^\nu.
    \end{equation}
    If $\ell=2$, then $\#\mu'=1$.
    But then since $\nu\ne(n-\mu_\ell)$, $\nu\not\!\unrhd\mu'$.
    Thus, $K_{\nu,\mu'}=0$, so $\mu_\ell K_{\nu,\mu'}\leq f^\nu$ holds either way.
    Therefore, by \Cref{kostka recursion},
    \begin{align*}
        (n-1)K_{\lambda,\mu} &= (n-1)\sum_\nu K_{\nu,\mu'}\\
        &= \sum_\nu ((n-\mu_\ell-1)K_{\nu,\mu'}+ \mu_\ell K_{\nu,\mu'}) \\
        &\leq \sum_\nu ((\ell-2)f^\nu +\mu_\ell K_{\nu,\mu'}) &\text{(by (\ref{eqn:IH}))} \\
        &\leq \sum_\nu ((\ell-2)f^\nu + f^\nu) &\text{(by (\ref{eqn:secondTerm}))} \\
        &\leq (\ell-1) f^\lambda. &\text{(by \Cref{dimension inequality})}
    \end{align*}
\end{proof}

\begin{prop}\label{main inequality}
	Let $n,d$ be nonnegative integers with $n\geq2$.
    Let $\lambda\vdash n$ with $\lambda\not=(n)$.
	If $n>d$, then
    $$\frac{F^{(n-1,1)}_d}{f^{(n-1,1)}}\geq\frac{F^\lambda_d}{f^\lambda}.$$
\end{prop}
\begin{proof}
	Recall that $f^{(n-1,1)}=n-1$ and $K_{(n-1,1),m(\ba)}=\#m(\ba)-1$ by \Cref{kostka props}(4) and \Cref{n-1 kostka}.
	Thus, proving the inequality amounts to showing that
	$$\sum_{\ba\vdash d,\#\ba=n}(\#m(\ba)-1)f^\lambda\geq \sum_{\ba\vdash d,\#\ba=n}(n-1)K_{\lambda,m(\ba)}.$$
	If $\lambda\ne(2,2)$, apply \Cref{kostka inequality} to conclude that the inequality holds for each summand, and thus for the sum.
	If $\lambda=(2,2)$ and $d=3$, then $m(\ba)\not=(2,2)$ for any $\ba\vdash 3$, so apply \Cref{kostka inequality} as before.
	If $\lambda=(2,2)$ and $d\le2$, a direct computation verifies the inequality.
\end{proof}

\begin{thm}\label{main theorem}
    Fix $n>d\ge 2$ and $1\leq r\leq P(d)$.
    Let $M$ be a maximal $r$-generated submodule of $R_d$.
    Then $\dim_k M\le\dim_k R_d-(P(d)-r)$ with equality if and only if 
   	\begin{equation}\label{var bound}
   	n\geq 1+\frac{\sum_{i=0}^{d-1}P(i)}{r}.
   	\end{equation}
\end{thm}
\begin{proof}
By \Cref{decomp of R_d},
	$$\dim_k R_d - \dim_k M= \sum_{\lambda\vdash n}(F^\lambda_d-\min\{F^\lambda_d,rf^\lambda\})f^\lambda.$$
	For $\lambda=(n)$, the summand is $P(d)-r$ by \Cref{kostka props} because any partition $\ba$ of $d$ can appear in the sum defining $F^{(n)}_d$ since $n\ge d$.
	For $\lambda\ne(n)$, the summand is $\ge0$ with equality precisely when $F^\lambda_d\le rf^\lambda$. 
	By \Cref{kostka partition sum} and \Cref{main inequality}, this occurs for all $\lambda\ne(n)$ exactly when $n\ge 1+\frac{\sum_{i=0}^{d-1}P(i)}{r}$.
\end{proof}

By \Cref{r gen characterization} and \Cref{decomp of R_d}, $\left\lceil\frac{F^\lambda_d}{f^{\lambda}}\right\rceil$ is the minimum number of $k\S_n$-module generators of the $\lambda$-isotypic component of $R_d$.
\Cref{main inequality} identifies the $(n-1,1)$-isotypic component as requiring the most orbit generators (excluding the $(n)$-isotypic component, which needs more generators).
Computational evidence suggests this is part of a more general phenomenon: the number of orbit generators appears to be monotone with respect to dominance order.

\begin{conj}
	Let $n,d$ be integers with $n,d\ge 2$.
	Let $\lambda,\mu\vdash n$.  If $\lambda\unrhd\mu$, then
	$$\left\lceil\frac{\sum_{\ba\vdash d,\#\ba=n}K_{\lambda,m(\ba)}}{f^{\lambda}}\right\rceil\ge\left\lceil\frac{\sum_{\ba\vdash d,\#\ba=n}K_{\mu,m(\ba)}}{f^{\mu}}\right\rceil$$
\end{conj}

\section{Applications to general symmetric ideals}\label{gsi section}

The following definition was introduced in \cite{SS} to model symmetric ideals generated by a collection of $r$ orbits of general polynomials of degree $d$.

\begin{defn}\label{general symmetric ideal}
	Let $\Gr_k(r,R_d)$ be the Grassmannian of $r$-dimensional linear subspaces of $R_d$.
	Define $\Phi:\Gr_k(r,R_d)\to \{\text{$r$-generated symmetric ideals}\}$ by
	$\Phi(L)=(L)_{\S_n}$.
	We say that a property $\mathscr{P}$ holds for a {\bf general $(r,d)$-symmetric ideal} if there exists a nonempty Zariski open set $U$ of $\Gr_k(r,R_d)$ such that for all $L\in U$, property $\mathscr{P}$ holds for $\Phi(L)=(L)_{\S_n}$. 
	
	If $r=1$ and thus $L\in U$ is spanned by a single polynomial $f$, the ideal $(f)_{\S_n}$ is called a {\bf general principal symmetric ideal}.
\end{defn}

\subsection{Relationship to maximal $r$-generated submodules}
General symmetric ideals can be difficult to work with as one needs to construct an explicit symmetric ideal to show that the Zariski open set is nonempty.
Fortunately, their properties are limited by the properties of maximal $r$-generated submodules of $R_d$.

\begin{prop} \label{general r gen is max}
    Let $I$ be a general $(r,d)$-symmetric ideal and let $M$ be a maximal $r$-generated submodule of $R_d$.
    Then $\dim_k I_d=\dim_k M$. In particular, $I_d$ is a maximal $r$-generated submodule of $R_d$.
\end{prop}
\begin{proof}
	Let $d'=\dim_k M$.
	Let
	$$U=\{L\mid \dim_k \Phi(L)_d=d'\}\subseteq\Gr_k(r,R_d).$$
	We first show that $U=\{L\in\Gr_k(r,R_d)\mid \dim_k \Phi(L)_d\geq d'\}$.
	One containment is obvious.
	For the other containment, notice that if $L\in\Gr_k(r,R_d)$, then $\Phi(L)$ is an $r$-generated symmetric ideal, so $\Phi(L)_d$ is an $r$-generated submodule and thus by \Cref{maximal r gen sub decomp} is isomorphic to a submodule of $M$.
	Hence, $\dim_k \Phi(L)_d\le \dim_k M=d'$.
	
	We now show that $U=\{L\mid \dim_k \Phi(L)_d\geq d'\}$ is a nonempty Zariski open set in $\Gr_k(r,R_d)$.
	Since $M$ is a maximal $r$-generated submodule of $R_d$, there is an $L\in \Gr_k(r,R_d)$ such that $M=\Phi(L)_d$.
	Thus, $\dim_k \Phi(L)_d=d'$, so $L\in U$.
	Let $N=\dim_k R_d$.
	Recall that $\Gr_k(r,R_d)$ is regarded as a subvariety of $\P^{\binom{N}{r}-1}$ via the Pl\"ucker embedding, which maps a subspace $L$ to the point in $\P^{\binom{N}{r}-1}$ whose coordinates are the $r\times r$ minors of a matrix whose row space is $L$.
	Throughout, we use the degree $d$ monomials, denoted $m_1,\dots,m_N$, as a basis over which to interpret the row space of a matrix.
	
	With $j_1<j_2<\cdots<j_r$, let $U_{j_1,\dots,j_r}\subseteq\P^{\binom{N}{r}-1}$ be the open set of points whose coordinate corresponding to the minor formed from columns $j_1,\dots,j_r$ is nonzero.
	It suffices to show that $U_{j_1,\dots,j_r}\cap U$ is open in the Grassmannian.
	If $L\in U_{j_1,\dots,j_r}$ and $C$ is a matrix whose row space is $L$, then by scaling the Pl\"ucker coordinates of $L$, one may assume the submatrix consisting of columns $j_1,\dots,j_r$ in $C$ is the identity matrix.
	It follows that the entries of $C$ in the remaining columns are Pl\"ucker coordinates of $L$.
	
	Let $f_1,\dots,f_r$ be a basis for $L$.
	Construct a matrix $D$ from $C$ whose rows correspond to the list of polynomials $\sigma\cdot f_i$ for each $\sigma\in \S_n$ and $1\leq i\leq r$, written in the fixed monomial basis.
	The entries of $D$ are the entries of $C$ permuted in various ways, and so the $d'\times d'$ minors of $D$ are polynomials in the entries of $C$.
	By homogenizing these polynomials with respect to the Pl\"ucker coordinate corresponding to $j_1,\dots,j_r$, we obtain homogeneous polynomials in the Pl\"ucker coordinates, say $g_1,\dots,g_t$, such that $U_{j_1,\dots,j_r}\cap U=U_{j_1,\dots,j_r}\backslash V(g_1,\dots,g_t)$.
\end{proof}

We now use \Cref{main inequality} and \Cref{r gen characterization} to provide an explicit construction of an $r$-generated symmetric ideal with maximal possible dimension in its minimum degree when \eqref{var bound} holds.

For each $\ba\vdash d$ with $\#\ba=n$, let $m_\ba$ be the sum of all monomials whose exponent vector is a permutation of $\ba$.  
These are the {\bf monomial symmetric polynomials}.  
For example, over $k[x_1,x_2,x_3]$, $m_{(3)}=x_1^3+x_2^3+x_3^3$.

\begin{constr}\label{r gen construction}
Let $A$ be a set of $r$ partitions of $d$ with $n$ parts.
Let $$J=(\bx^\ba-\s\cdot\bx^\ba : \s\in\S_n, \ba\vdash d,\#\ba=n,\ba\notin A)+(\bx^\ba : \ba \in A)_{\S_n}$$ 
$$W=\langle m_\ba : \ba\vdash d,\#\ba=n,\ba\notin A\rangle_{k\S_n}=\bigoplus_{\ba\vdash d,\#\ba=n,\ba\notin A}\langle m_\ba\rangle_{k\S_n}$$
\end{constr}
\begin{lem}\label{W structure}
	Let $W$ be as in \Cref{r gen construction}.
    If $n\geq d$, then $W\cong (P(d)-r)S^{(n)}$ and so $\dim_k W=P(d)-r$.
\end{lem}
\begin{proof}
    Since $m_\ba$ is symmetric, $\langle m_\ba\rangle_{k\S_n}$ is the trivial representation of $\S_n$, so $\langle m_\ba\rangle_{k\S_n}\cong S^{(n)}$.
    Since $n\geq d$, every partition of $d$ not in $A$ appears in the direct sum defining $W$, so $W\cong (P(d)-r)S^{(n)}$.
    Since $\dim_k S^{(n)}=1$, $\dim_k W=P(d)-r$.
\end{proof}
\begin{lem}\label{RJW}
	Let $J,W$ be as in \Cref{r gen construction}.
	Then $R_d\cong J_d\oplus W$ as $k\S_n$-modules.
\end{lem}
\begin{proof}
	We first show that $R_d\subseteq J_d+W$.
	Let $\bx^\ba$ be a monomial in $R_d$.
	Let $C$ be a complete set of representatives for the quotient of $\S_n$ by the stabilizer of $\ba$.
	That is, $C\subseteq\S_n$ is such that $m_\ba=\sum_{\s\in C}\s\cdot\bx^\ba$.
	Then
	$$m_\ba+\sum_{\s\in C}(\bx^\ba-\s\cdot\bx^\ba)=\sum_{\s\in C}\bx^\ba=|C|\bx^\ba$$
	is an element of $J_d+W$.
	Since $|C| \mid n!$, $|C|$ is invertible in $k$, so $\bx^\ba\in J_d+W$, as desired.
	
	The sum of the coefficients of the monomials in each $m_\ba$ is the cardinality of the orbit of $\ba$ under the action of $\S_n$ and thus it divides $n!$ and so is nonzero in $k$.
	Furthermore, for any $f\in W$, the sum of the coefficients of monomials in $f$ whose exponent vector is a permutation of $\ba\notin A$ is either nonzero or no such monomials appear in $f$.
	However, if $f\in J_d$, then this same sum must be 0, since this property holds for the generators of $J_d$.
	Thus, if $f\in J_d\cap W$, then the only monomials which can appear in $f$ are a permutation of $\bx^\ba$ for some $\ba\in A$, but no such monomial appears in any element of $W$, so $f=0$.
\end{proof}

\begin{prop}\label{J props}
Let $J$ be as in \Cref{r gen construction} and let $n\ge d$.
Then $J$ is a symmetric ideal generated in degree $d$ such that $\dim_k J_d=\dim_k R_d-(P(d)-r)$.
Furthermore, if $n>d$ and $n\geq 1+\frac{\sum_{i=0}^{d-1}P(i)}{r}$, then $J$ is $r$-generated.
\end{prop}
\begin{proof}
	The first statement is immediate from the definition of $J$, \Cref{W structure} and \Cref{RJW}.
	
	For the second statement, observe that by \Cref{decomp of R_d}, \Cref{W structure}, and \Cref{RJW}, $J_d$ decomposes as 
	$$rS^{(n)}\oplus\bigoplus_{\lambda\vdash n, \lambda\ne(n)}F^\lambda_d S^\lambda$$
	By \Cref{kostka partition sum} and \Cref{main inequality}, for $\lambda\not=(n)$,
    $$r\ge\frac{\sum_{i=0}^{d-1}P(i)}{n-1}=\frac{F^{(n-1,1)}_d}{f^{(n-1,1)}}\geq \frac{F^\lambda_d}{f^\lambda}.$$
	So $J$ is an $r$-generated symmetric ideal by \Cref{r gen characterization}.
\end{proof}

\subsection{Effective results}

In \cite{HSS}, the authors give the Hilbert function, Betti table, and graded minimal free resolution of the quotient of $R$ by a general principal symmetric ideal in a sufficiently large number of variables.
When $n$ is large enough, the authors construct an ideal with certain properties to show that the Zariski open set in which their conclusion is true is nonempty.
However, they did not provide an explicit bound on the number of variables needed for their theorem.
In this section, we use \Cref{main inequality} and \Cref{r gen characterization} to make their theorem effective by giving a sharp bound on the number of variables needed for their conclusion to hold.

A careful analysis of the proof of \cite[Theorem 8.4]{HSS} reveals that for the proof to go through, we need an ideal $J$ with the following three properties:
\begin{enumerate}
    \item $J$ is a principal symmetric ideal
    \item $\dim_k J_d=\dim_k R_d-(P(d)-1)$
    \item If $J=(f)_{\S_n}$, then for each $\ba\vdash d$,$\#\ba=n$, $\ba\not=(d,0,\dots,0)$, the sum of the coefficients of monomials in $f$ whose exponent vector is a permutation of $\ba$ is 0, and the sum of the coefficients of monomials in $f$ whose exponent vector is a permutation of $(d,0,\dots,0)$ is nonzero.
\end{enumerate}

We have shown in \Cref{J props} the existence of an ideal with the first two properties.
We now establish that it satisfies the last as well when $r=1$.

\begin{prop}\label{J technical requirement}
	Let $J$ be as in \Cref{r gen construction} with $A=\{(d,0,\dots,0)\}$.
	Then (3) holds.
\end{prop}
\begin{proof}
	Suppose that $J=(f)_{\S_n}$.
	Then $f\in J_d$ and in the proof of \Cref{RJW}, we established that every element of $J_d$ satisfies the first part of (3).
	By similar reasoning, if the sum of the coefficients of monomials in $f$ whose exponent vector is a permutation of $(d,0,\dots,0)$ is 0, then this will be true for every element of $\langle f\rangle_{k\S_n}=J_d$.
	But $x_1^d\in J_d$, so this does not happen.
\end{proof}

One change that we make to the theorem in \cite{HSS} is to allow for char$(k)>n$.
This can be done because they use the assumption that char$(k)=0$ in only two places.
First, in the construction of their symmetric ideal, which we are not using.
Second, in \cite[Proposition 7.15]{HSS}.
However, a careful analysis shows that the proof they give also holds when char$(k)>n$.

\begin{thm}[Strengthening of {\cite[Theorem 8.4]{HSS}}]\label{8.4 strengthening}
    Assume $k$ is infinite with either char$(k)=0$ or char$(k)>n$.
    Fix $d\geq 2$ and $n>d$.
    Let $I\subseteq k[x_1,\dots,x_n]$ be a general principal symmetric ideal generated in degree $d$.
    Then $I$ has the following properties if and only if $n\geq 1+\sum_{i=0}^{d-1}P(i)$:
    \begin{enumerate}
        \item the Hilbert function of $A=R/I$ is given by 
        $$H_A(i):=\dim_k A_i=\left\{\begin{array}{ll}
            \dim_k R_i & \text{if }i\leq d-1 \\
             P(d)-1 & \text{if }i=d \\
             0 & \text{if }i>d
        \end{array}\right.$$

        \item the Betti table of $A$ has the form given in \cite[Theorem 8.4]{HSS}
        \item the graded minimal free resolution of $A$ has $\S_n$-equivariant structure described by the $k\S_n$-irreducible decompositions for the modules $\Tor^R_i(A,k)$ given in \cite[Theorem 8.3]{HSS}.
    \end{enumerate}
    Moreover, the Poincar\'{e} series of all finitely generated graded $A$-modules are rational, sharing a common denominator.
    When $d>2$, $A$ is Golod.
    When $d=2$, $A$ is Gorenstein and Koszul.
\end{thm}
\begin{proof}
    If $n\geq1+\sum_{i=0}^{d-1}P(i)$, then by \Cref{J props} and \Cref{J technical requirement}, $J$ has the required properties, so we can replicate the proof given in \cite{HSS} starting with their Theorem 8.1, replacing their constructed polynomial with a generator of $J$.
    Conversely, if $I$ has the properties above, then $\dim_k I_d=\dim_k R_d-(P(d)-1)$.
    Thus, by \Cref{main theorem} and \Cref{general r gen is max}, $n\geq 1+\sum_{i=0}^{d-1}P(i)$.
\end{proof}

In \cite{HSS}, they analyzed the case $d=3$ in further detail and determined that a tight bound for the conclusion holding is $n\geq 5$.
This analysis agrees with \Cref{8.4 strengthening}.
The recent paper \cite{SS} generalizes the results of \cite{HSS} to general $(r,d)$-symmetric ideals and uses the results of this paper to give an analogous effective bound on the number of variables.

\appendix

\section{Calculations for the Kostka numbers inequality}
\label{kostka inequality appendix}

Here we give the details to complete the proof of Case 3 and Case 4 in \Cref{kostka inequality}.
The following lemma finishes the proof of Case 4.
\begin{lem}\label{inequality proof}
	Let $n\geq5$ and $\ell,p\geq 2$ such that $n\geq \ell\cdot p$.
	Then
	$$\frac{(n-1)(p+\ell-2)!}{p!(\ell-2)!}\leq\frac{(\ell-1)n!(n-2p+1)}{p!(n-p+1)!}.$$
\end{lem}
\begin{proof}
	We consider the three cases $\ell\geq 3$ and $p\geq 3$,  $\ell=2$ and $p\geq3$ and $\ell\geq 2$ and $p=2$.

    {\em Case 1:} Suppose that $\ell\geq 3$ and $p\geq 3$.
    By assumption,
    $$n\geq \ell\cdot p=(\ell-2)(p-1)+2p+\ell-2.$$
    Since $\ell\geq 3$ and $p\geq 3$, $(\ell-2)(p-1)\geq 1$, so $n\geq 2p+\ell-1$, thus $n-2p+1\geq \ell$.
    Therefore,
	\begin{align*}
        \frac{(n-1)(p+\ell-2)!}{p!(\ell-2)!} 
        &= \frac{n-1}{p!}(p+\ell-2)(p+\ell-3)\cdots \ell(\ell-1) \\
        &\leq \frac{(\ell-1)(n-1)(n-2p+1)}{p!}(p+\ell-2)(p+\ell-3)\cdots(\ell+1) \\
        &\leq \frac{(\ell-1)(n-1)(n-2p+1)}{p!}n(n-2)\cdots(n-p+2) \\
        &= \frac{(\ell-1)n!(n-2p+1)}{p!(n-p+1)!}.
    \end{align*}    
	The second inequality compares factors, using $p+\ell-2\le n$ and $p+\ell-3-i\le n-2-i$ for $0\le i \le p-4$.

    {\em Case 2:} Suppose that $\ell=2$ and $p\geq 3$.
    Then the inequality reduces to 
    $$n-1\leq\frac{n!(n-2p+1)}{p!(n-p+1)!}.$$
    We start from the right hand side this time.
    Note that $\binom{n}{p}\geq\binom{n}{3}$ as $p\geq 3$ and $p\leq\frac{n}{2}$.
    Thus, 
    \begin{align*}
        \frac{n!(n-2p+1)}{p!(n-p+1)!} 
        &= \binom{n}{p}\frac{n-2p+1}{n-p+1} \\
        &\geq \binom{n}{3}\frac{n-2p+1}{n-p+1} \\
        &= \frac{n(n-1)(n-2)(n-2p+1)}{6(n-p+1)} \\
        &\geq \frac{n(n-1)(n-2p+1)}{6}.
    \end{align*}
    Since $p\geq 3$ and $\ell=2$, $n\geq \ell\cdot p\geq6$.
    Hence, 
    $$\frac{n!(n-2p+1)}{p!(n-p+1)!}\geq\frac{n(n-1)(n-2p+1)}{6}\geq (n-1)(n-2p+1)\geq n-1.$$

    {\em Case 3:} Suppose that $\ell\geq2$ and $p=2$.
    The inequality becomes 
    $$\frac{(n-1)\ell!}{2(\ell-2)!}\leq \frac{(\ell-1)n!(n-3)}{2(n-1)!}.$$
    Recall that $n\geq \ell\cdot p=2\ell$ and $n\ge 5$.
    Hence, $n-3\ge \ell$.
	Therefore,
    \begin{align*}
		\frac{(\ell-1)n!(n-3)}{2(n-1)!}
		&\geq \frac{(\ell-1)(n-1)\ell}{2} \\
		&= \frac{(n-1)\ell!}{2(\ell-2)!}. \qedhere
    \end{align*}
\end{proof}

Here we complete the computations for Case 3 of \Cref{kostka inequality} by verifying that \eqref{kostka term inequality} holds when $\mu=(2,2,1)$ and $\mu=(2,2,2)$ for any choice of $\lambda\not=(n)$, where $\mu\vdash n$.
Several choices of $\lambda$ are covered by Case 1 of \Cref{kostka inequality}, and $\lambda=(4,2),\mu=(2,2,2)$ is covered by Case 4.
The following table shows the values for the remaining cases.\\

\begin{center}
\begin{tabular}{l|c|c|c|c|c}
    $\lambda$ & $\mu$ & $K_{\lambda,\mu}$ & $f^\lambda$ & $(n-1)K_{\lambda,\mu}$ & $(\#\mu-1)f^\lambda$ \\\hline  
    (3, 2)    & (2, 2, 1) & 2 & 5 & 8 & 10 \\
    (3, 1, 1) & (2, 2, 1) & 1 & 6 & 4 & 12 \\
    (2, 2, 1) & (2, 2, 1) & 1 & 5 & 4 & 10 \\
    (4, 1, 1) & (2, 2, 2) & 1 & 10 & 5 & 20 \\
    (3, 3)    & (2, 2, 2) & 1 & 5 & 5 & 10 \\
    (3, 2, 1) & (2, 2, 2) & 2 & 16 & 10 & 32 \\
    (2, 2, 2) & (2, 2, 2) & 1 & 5 & 5 & 10 \\
\end{tabular}
\end{center}

\phantom{.}\\
The values can be checked by hand or by computer \cite{GithubRepo}.

\begin{rem}
A complete formalization of the proof of \Cref{main inequality} presented in Section 4 has been written by the author in Lean using Mathlib.
This formalization effort actually uncovered an error in the original proof of \Cref{inequality proof}, which has now been corrected.
The code is located in \href{https://github.com/NoahW314/KostkaNumbers}{this GitHub repository}.
\end{rem}

\noindent
{\bf Acknowledgments.}
We thank Alexandra Seceleanu for her advice and helpful comments, which greatly improved this manuscript.
We also thank the anonymous reviewers for their useful feedback.
We acknowledge Macaulay2, which was used to perform many computations in the earlier stages of this research, and the OEIS for its help in identifying A000070.

\printbibliography

\end{document}